\documentclass[sn-aps]{sn-jnl}
\usepackage{graphicx}

\graphicspath{ {./} }
\usepackage{epsfig} 
\usepackage{svg}

\usepackage{amsmath}
\usepackage{amssymb}
\usepackage{graphicx}
\usepackage{color,soul}
\usepackage{amsthm}
\usepackage{algpseudocode}
\usepackage{algorithm}
\usepackage{verbatim}

\usepackage{caption}
\usepackage{subcaption}
\usepackage{setspace}
\usepackage{graphicx}%
\usepackage{multirow}%
\usepackage{amsmath,amssymb,amsfonts}%
\usepackage{amsthm}%
\usepackage{mathrsfs}%
\usepackage[title]{appendix}%
\usepackage{xcolor}%
\usepackage{textcomp}%
\usepackage{manyfoot}%
\usepackage{booktabs}%
\usepackage{algorithm}%
\usepackage{algorithmicx}%
\usepackage{algpseudocode}%
\usepackage{listings}%
\usepackage{lineno}
\usepackage{caption}

\title{Computing a Sparse Approximate Inverse on Quantum Annealing Machines }

\author[1]{\fnm{Sanjay} \sur{Suresh}}\email{ssuresh27@wisc.edu}

\author*[2]{\fnm{Krishnan} \sur{Suresh}}\email{ksuresh@wisc.edu}

\affil[1]{\orgdiv{Computer Science}, \orgname{University of Wisconsin, Madison}, \orgaddress{\street{1210 W. Dayton Street}, \city{Madison}, \postcode{53706}, \state{WI}, \country{USA}}}

\affil*[2]{\orgdiv{Mechanical Engineering}, \orgname{University of Wisconsin, Madison}, \orgaddress{\street{1513 University Avenue}, \city{Madison}, \postcode{53706}, \state{WI}, \country{USA}}}

\begin{document}

\abstract
{ Many engineering problems involve solving large linear systems of equations. Conjugate gradient (CG) is one of the most popular iterative methods for solving such systems. However, CG typically requires a good preconditioner to speed up convergence. One such preconditioner is the sparse approximate inverse (SPAI). 

In this paper, we explore the computation of an SPAI on quantum annealing machines by solving a series of quadratic unconstrained binary optimization (QUBO) problems. Numerical experiments are conducted using both well-conditioned and poorly-conditioned linear systems arising from a 2D finite difference formulation of the Poisson problem. 
}

\maketitle

\section{Introduction}

Many engineering problems result in linear systems of equations \cite {zienkiewicz2005finite} of the form:
\begin{equation}\label{bigMatrixProblem}
 \mathbf K  \mathbf u = \mathbf f
\end{equation}
where $ \mathbf K  $ is a symmetric, positive-definite, sparse $N \times N$ matrix, $  \mathbf u$ is the unknown field (for example, the temperature field), and $  \mathbf f$ is the applied force (for example, the heat-flux). Solving such linear systems for large $N$ is a computationally intensive task \cite {gould2007numerical}, and can be time-consuming on classical computers. Quantum computers have been proposed as an alternate since they can potentially accelerate the computation; see \cite{tosti2022review, wang2023opportunities} for recent reviews on the potential role of quantum computers in engineering. 

In particular, the Harrow-Hassidim-Lloyd (HHL) algorithm is a landmark strategy for solving linear systems of equations on \emph {quantum-gate computers}. In theory, it offers an exponential speed-up over classical algorithms \cite{harrow2009quantum}, and it has been further improved recently \cite{ambainis2010variable,childs2017quantum, liu2022survey}.  However, due to the accumulation of errors in current noisy intermediate-scale quantum (NISQ) computers  \cite{preskill2018quantum}, the HHL algorithm and its variants are limited, in practice, to very small ($ N < 5$) systems. Other strategies such as Grover's algorithm \cite{srinivasan2017solving} and quantum approximate optimization algorithm (QAOA)  \cite{an2022quantum} have been proposed to solve linear systems on quantum gate computers, but they suffer from the same limitation. Hybrid solvers such as the variational quantum linear solver \cite{bravo2019variational} have also been proposed to mitigate some of these challenges.

In parallel, \emph {quantum annealing machines}, such as the D-wave systems with several thousand qubits \cite{shin2014quantum}, have also been proposed for solving linear systems since they are less susceptible to noise \cite{hauke2020perspectives, yarkoni2022quantum}. The basic principle is to pose the solution of a linear system of equations as a minimization problem and then convert this into a series of quadratic unconstrained binary optimization (QUBO) problems. For example, O'Malley and Vesselinov used a least-squares formulation and a finite-precision qubit representation to pose QUBO problems \cite{o2018nonnegative}.  Borle and Lomonaco carried out a theoretical analysis of this approach \cite{borle2019analyzing}. Park et. al. showed how QUBO problems can be simplified using matrix congruence \cite{park2021application}, while its application in solving 1D Poisson problems was demonstrated in \cite{conley2023quantum}. If the matrix is positive-definite, which is often the case in engineering problems, it is much more efficient to use a potential-energy formulation, rather than the least-squares formulation, to pose QUBO problems. Using the potential-energy formulation, Srivastava et. al. described a box algorithm to solve QUBO problems arising from finite element analysis of one-dimensional differential equations \cite{srivastava2019box}.

However, despite these advances,  only small-size ($ N < 100 $) problems have been demonstrated on current quantum annealing machines, while most engineering problems of practical interest result in much larger-size ($ N > 10^5 $) matrices. To address this gap, we propose here an alternate paradigm where \emph {quantum annealing machines are used, not to directly solve such large linear systems, but to compute sparse preconditioners (SPAI)}. Once again, the strategy is to rely on the potential energy formulation to pose a series of quadratic unconstrained binary optimization (QUBO) problems. However, the sparsity can be fully exploited to significantly reduce the size of the QUBO problems, making it more amenable to quantum computing. The computed SPAI can then be used as a preconditioner to rapidly solve large systems of equations using the conjugate gradient (CG) method on classical computers.  The efficacy of this approach is demonstrated by solving both well-conditioned and poorly-conditioned linear systems of equations arising from finite difference formulation of 2D Poisson problems. Furthermore, by exploiting the well-structured nature of the finite-difference formulation, we show how one can compute the SPAI in constant time, independent of the size ($N$) of the linear system. 

\section{Proposed Methodology} \label{sec:overviewQCalgorithms}

\subsection{Poisson Problem}

To provide a context to this paper, we consider solving a Poisson problem in 2D, governed by \cite{langtangen2017finite}:
\begin{equation} 
\label{Poisson} 
k \bigg(\frac{\partial^2 u}{\partial x^2} + \frac{\partial^2 u}{\partial y^2} \bigg)  = -f
\end{equation} 
where $u$ is the unknown (temperature) field, $f$ is the (heat) source, and $k$ is the underlying material property (example: conductivity). We will assume that the field takes a value of zero on the boundary. 

A classic approach for solving Eq. \ref{Poisson} is the finite difference method \cite{langtangen2017finite}, where the geometry is sampled by a uniform grid, over which the field $u$ is to be determined. Then, the partial derivatives are approximated as follows \cite{langtangen2017finite}:

\begin{equation}
    \frac{\partial^2 u}{\partial x^2} \Big\rvert_{m,n} \approx \frac{u_{m+1,n}-2u_{m,n}+u_{m-1,n}}{\Delta x^2}
\end{equation}
\begin{equation}
    \frac{\partial^2 u}{\partial y^2} \Big\rvert_{m,n} \approx \frac{u_{m,n+1}-2u_{m,n}+u_{m,n-1}}{\Delta y^2}
\end{equation}
Substituting the above approximations in Eq. \ref{Poisson}, and with $\Delta x = \Delta y = h$, results in:
\begin{equation}
\label{fdPoisson} 
4ku_{m,n} -  k\left( u_{m-1,n} + u_{m+1,n} + u_{m,n-1} + u_{m,n+1}\right) = f_{m,n}h^2
\end{equation}
When applied at all grid points, and after eliminating the rows and columns corresponding to the boundary nodes, we arrive at the linear system in Eq. \eqref{bigMatrixProblem}. As one can observe, the resulting $\mathbf K$ matrix is symmetric. Further, for this problem, the sparsity is $5$, i.e., $\mathbf K$ has at most 5 entries in any row or column.  Once the rows and columns corresponding to the boundary nodes are eliminated, and the resulting $\mathbf K$ matrix can be shown to be positive-definite \cite{langtangen2017finite}.

\begin{figure}[H]
 	\begin{center}
	\includegraphics[width=0.4\textwidth]{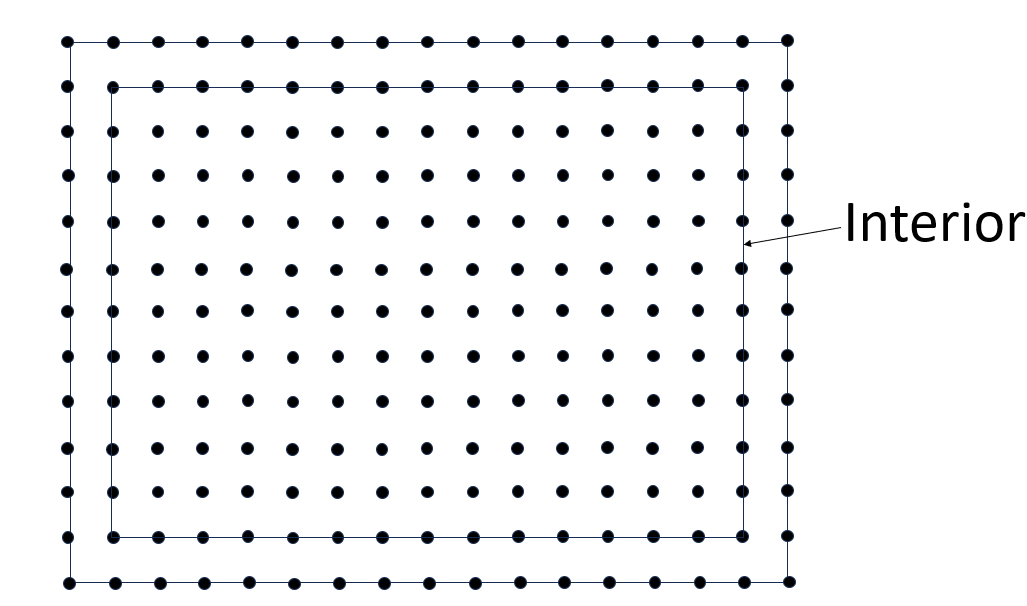}
    \caption{Finite-difference grid over a rectangle with zero Dirichlet boundary conditions.}
	\label{fig:fdGrid}
	\end{center}
 \end{figure}

While the finite difference method applies to arbitrary domains, we will focus on rectangular domains for simplicity; see Fig \ref{fig:fdGrid}. Since the field is assumed to take a value of zero on the boundary, we will only solve for the field in the interior, i.e., $N$ denotes the number of interior nodes.

Finally, if the material property ($k$) is not a constant over the domain,  Eq. \ref{fdPoisson} can be generalized to
\begin{equation}
\label{fdPoissonGeneral} 
4k_{m,n} u_{m,n} -  \left( k_{m,n}^{m-1,n} u_{m-1,n} + k_{m,n}^{m+1,n}u_{m+1,n} + ...\right) = f_{m,n}h^2
\end{equation}
where  $k_{m,n}$ is the average $k$ over the (4 or less) elements surrounding the node $(m,n)$, while $k_{m,n}^{m-1,n}$  is the average $k$ over the (2 or 1) elements adjacent to the edge joining node $(m,n)$ and $ (m-1,n)$. The resulting $\mathbf  K $ matrix is still symmetric, sparse, and positive definite, but can become poorly conditioned depending on the material distribution (see Section \ref{sec:Experiments}).

\subsection{Classic Linear Solvers}

The two main strategies for solving Eq.\eqref{bigMatrixProblem} are: direct and iterative.  Direct solvers usually rely on LU-decomposition \cite{bollhofer2020state}, and generally require a significant amount of memory for large systems. Iterative solvers require less memory but converge to a solution gradually. The rate of convergence of iterative methods depends on several factors including the condition number of $\mathbf K$, sparsity, etc \cite{axelsson2007solution}. One of the most popular iterative methods is conjugate gradient whose run-time complexity is given by \cite{shewchuk1994introduction, nazareth2009conjugate}:
\begin{equation}\label{CGComplexity}
    CG \sim  O\left( N s \sqrt{\kappa} \log{ \frac 1 \epsilon_{cg}} \right)
\end{equation}
where $N$ is the dimension of $ \mathbf K  $, $s$ is the sparsity ($s = 5$ for the problem described above),  $\kappa$ is the condition number, and $\epsilon _{cg}$ is the desired residual error. Thus, for poorly conditioned systems (i.e. when $\kappa$  is large), CG does not perform well, and preconditioners are essential. Several preconditioners are widely used today \cite{chow2000priori, benzi2002preconditioning,wathen2015preconditioning}; these include Jacobi, incomplete Cholesky, sparse approximate inverse, etc. In this paper, we will rely on the sparse approximate inverse of $ \mathbf K  $, and we propose a simple algorithm to compute this preconditioner on quantum annealing computers using a QUBO formulation. 

\subsection{Computing a Sparse Approximate Inverse}
Note that $\mathbf K^{-1}$, the exact inverse of $\mathbf K$,  can be computed by solving 
\begin{equation}\label{exactInverse} 
\mathbf K \mathbf y_j = \mathbf e_j, j =0,1 ... N-1
\end{equation}
where $\mathbf e_j$ is the unit vector corresponding to the $j^{th}$ dimension. In general, $\mathbf K^{-1}$ will be dense. Our objective is to compute a \emph{sparse approximate inverse} (SPAI) $\mathbf M$.  There are various  \emph{a priori} and adaptive techniques for forcing sparsity on $\mathbf M$; see \cite{chow2000priori, benzi2002preconditioning}. We will use a well-known \emph{a priori} technique where the sparsity pattern of $\mathbf K$ is imposed on  $\mathbf M$ \cite{benzi2002preconditioning}.  In other words, to compute  $\mathbf M$, we once again  solve
\begin{equation}\label{approxInverse} 
\mathbf K \mathbf m_j = \mathbf e_j, j =0,1 ... N-1
\end{equation}
but with the constraint that  $\mathbf m_j$ must have the same sparsity pattern as the  $j^{th}$ column of $ \mathbf K  $. 

To compute $\mathbf m_j$, let the row-sparsity index of the $j^{th}$ column of $\mathbf K$ be $ \mathbf s$, i.e., $ \mathbf s[i]$ stores the $i^{th}$ non-zero row of $\mathbf K[:,j]$, $0 \le i \le s-1$, and $ s\ll N$. We can rearrange $\mathbf m_j$ as 
\begin{equation}
\mathbf m_j = \begin{bmatrix}
    \mathbf {\hat{m}}_j \\
    \mathbf 0
\end{bmatrix}
\end{equation}
where $\mathbf {\hat{m}}_j$ is of length $s$. We can also rearrange $\mathbf K$ using the same reordering as: 
\begin{equation}\label{K_extraction}
\mathbf K = \begin{bmatrix}
    \mathbf A_j & \mathbf B_j\\
    \mathbf B_j^T & \mathbf C_j 
\end{bmatrix}
\end{equation}
where $\mathbf A_j$ is  $s \times s$, and $\mathbf C_j$ is $(N-s) \times (N-s)$, and rearrange $\mathbf e_j$ as
\begin{equation}
\mathbf e_j = \begin{bmatrix}
    \mathbf {\hat{e}}_j \\
    \mathbf 0
\end{bmatrix}
\end{equation}
Thus  Eq. \ref{approxInverse}  reduces to:
\begin{equation}
\label{hamiltonian5} 
    \begin{bmatrix}
    \mathbf A_j & \mathbf B_j\\
    \mathbf B_j^T & \mathbf C_j 
\end{bmatrix} \begin{bmatrix}
    \mathbf {\hat{m}}_j \\
    \mathbf 0
\end{bmatrix} =  \begin{bmatrix}
    \mathbf {\hat{e}}_j \\
    \mathbf 0
\end{bmatrix}
\end{equation}
Discarding the last $(N-s)$ equations, we arrive at : 
\begin{equation}
\label{reducedProblem} 
\mathbf A_j \mathbf {\hat{m}}_j = \mathbf {\hat{e}}_j
\end{equation}
where $\mathbf A_j$ is a small $s \times s$ matrix constructed from the sparsity pattern of the $j^{th}$ column. Eq. \ref{reducedProblem} must be posed and solved for each column $j$ of $\mathbf K$. However, we show in Section \ref{sec:NodeMapping} that we only need to solve for a few columns by exploiting the structured nature of the finite difference formulation.

\subsection{QUBO Formulation}
Observe that solving Eq. \eqref{reducedProblem} is equivalent to minimizing the  potential energy (we have dropped the subscript $j$ to avoid clutter):
\begin{equation}
\label{hamiltonian2} 
    \min_{\mathbf {\hat{m}}} \Pi =  \frac 1 2 \mathbf {\hat{m}}^T\mathbf A \mathbf {\hat{m}} -  \mathbf {\hat{m}}^T \mathbf {\hat{e}} 
\end{equation}
This is a quadratic unconstrained minimization problem involving real variables $ \mathbf {\hat{m}}  $. To solve this on a quantum annealing machine, we represent each real component ${\hat{m}}_i$ using qubit variables. A well known strategy is the radix representation, also referred to here as the \emph {box representation} \cite{o2018nonnegative, rogers2020floating}: 
\begin{equation}\label{boxRepresentation}
\mathbf {\hat m} = \mathbf c + L(-2\mathbf q_1 + \mathbf q_2)
\end{equation}
where $\mathbf c $ and $L $ are real value parameters that we will iteratively improve, while $\mathbf q_1 $ and $\mathbf q_2 $ are qubit vectors of length $s$ each, i.e., a total of  $2s$ qubits is used to capture  $ \mathbf {\hat{m}} $.

Since $\mathbf {\hat m} $  is linear in $\mathbf q_1 $ and $\mathbf q_2 $, substituting Eq. \eqref{boxRepresentation} into Eq. \eqref{hamiltonian2} will lead to a quadratic unconstrained \emph {binary} optimization (QUBO) problem:
\begin{equation}
\label{hamiltonianQUBOInitial} 
    \min_{\mathbf q = \{ \mathbf q_1 , \mathbf q_2 \} } \Pi = \frac 1 2 \mathbf q^T \mathbf Q' \mathbf q  + \mathbf q^T \mathbf d
\end{equation}

Further, since the qubit variables can only take the binary values $0$ or $1$, the linear term can be absorbed into  the quadratic term \cite{zaman2021pyqubo}, resulting in the standard form: 
\begin{equation}
\label{hamiltonianQUBO} 
    \min_{\mathbf q = \{ \mathbf q_1 , \mathbf q_2 \} } \Pi = \frac 1 2 \mathbf q^T\mathbf Q \mathbf q 
\end{equation}
where $ \mathbf Q$ is symmetric (but not positive definite).  

The overall strategy therefore is as follows: For the first column of $\mathbf K $, the parameters  $\mathbf c$ and $L$ are initialized to $\mathbf c = \mathbf 0 $ and $L = 0$. The resulting QUBO problem in Eq. \ref{hamiltonianQUBO} is solved on a quantum annealing machine. Then $\mathbf c$ and $L$  in Eq. \ref{boxRepresentation} are updated via the \emph {sparse box algorithm}, discussed in Section \ref{sec:Algorithms}. The process is repeated until convergence is reached (in typically $25 \sim 35$ iterations; see Section \ref {sec:Experiments}). Then, the next column of $\mathbf K $ is processed until the entire $\mathbf M $ matrix is constructed.

\subsection{Node Mapping for a Structured Grid}
\label{sec:NodeMapping}
\begin{figure}[H]
 	\begin{center}
	\includegraphics[width=0.3\textwidth]{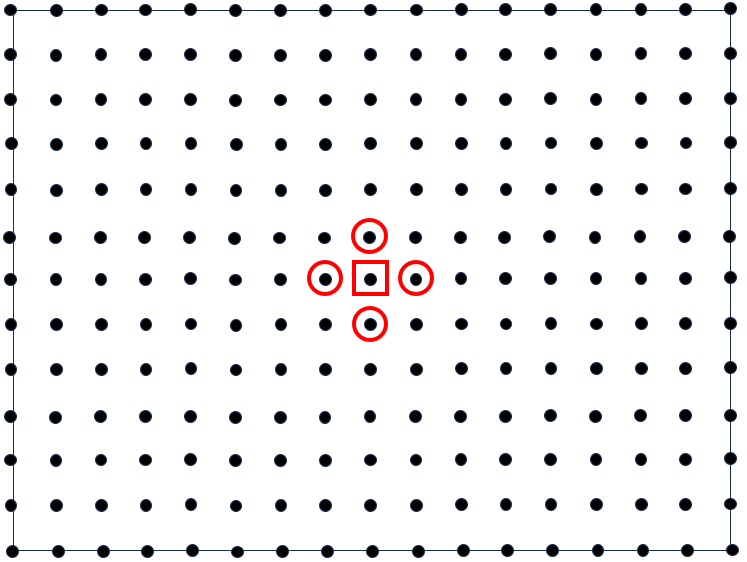}
    \caption{A  node and its 4 neighboring nodes.}
	\label{fig:fdGridNodeNeighbors}
	\end{center}
 \end{figure}

\begin{figure}[H]
 	\begin{center}
	\includegraphics[width=0.3\textwidth]{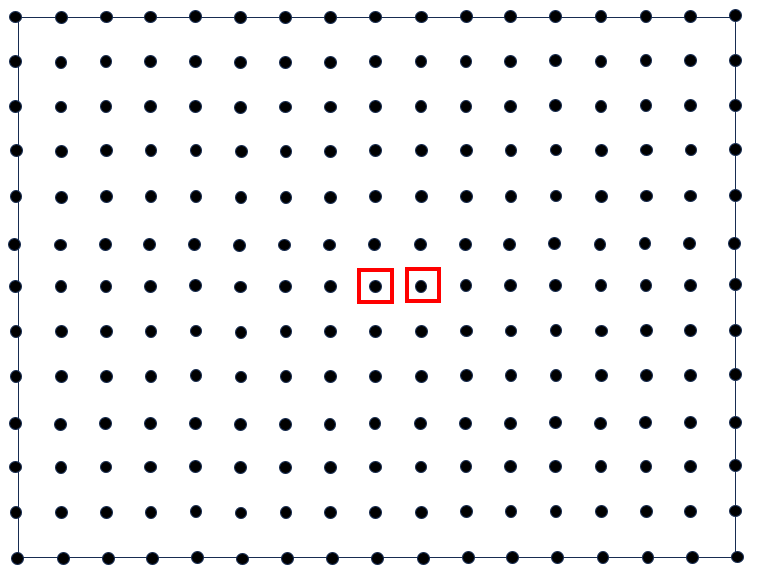}
    \caption{Two nodes with identical columns in $\mathbf M $.}
	\label{fig:fdGridNodeMapping}
	\end{center}
 \end{figure}

\begin{figure}[H]
 	\begin{center}
	\includegraphics[width=0.3\textwidth]{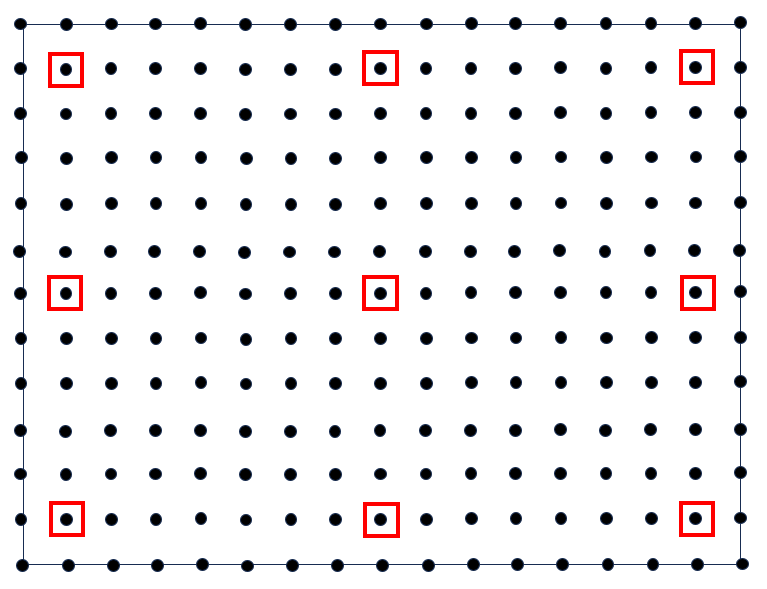}
    \caption{For a single material, only 9 nodes need to be considered.}
	\label{fig:fd9IndependentNodes}
	\end{center}
 \end{figure}

\begin{figure}[H]
 	\begin{center}
	\includegraphics[width=0.3\textwidth]{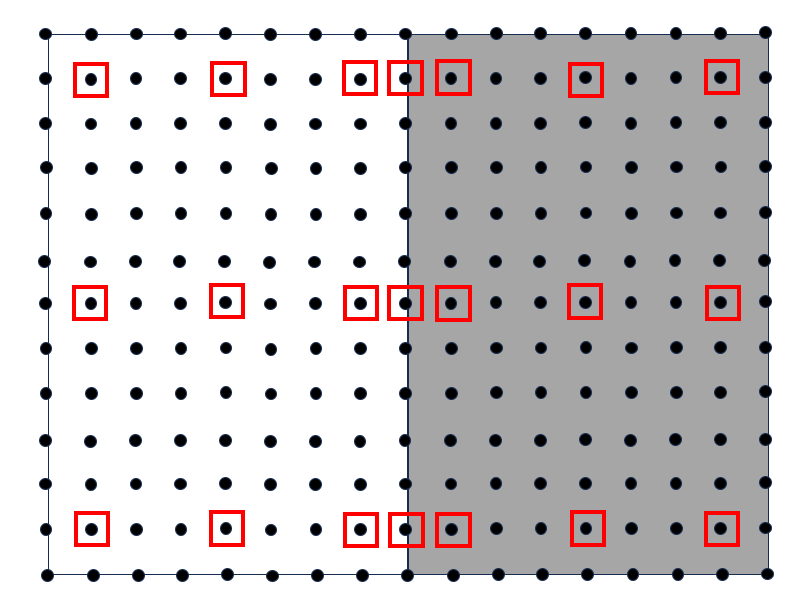}
    \caption{For this two-material configuration, only 21 nodes need to be considered.}
	\label{fig:2materials}
	\end{center}
 \end{figure}

For a generic $\mathbf K $ matrix, one must explicitly process each of the $N$ columns. However, for the structured finite difference grid, one can significantly reduce the computation. Observe that each column of  $\mathbf K $  corresponds to a  unique node in the grid. Now consider a typical node highlighted using a square box in Fig. \ref {fig:fdGridNodeNeighbors}. The  $5 \times 5$ matrix $\mathbf A $ in Eq. \ref{hamiltonian2} corresponding to this node (column), is entirely determined by the rows and columns of $\mathbf K $ associated with this node and the 4 neighboring nodes (highlighted using circles in Fig. \ref {fig:fdGridNodeNeighbors}). From Eq. \ref{fdPoissonGeneral}, we observe that the diagonal entries of $\mathbf A $ depend only on the material property ($ k_{m,n} $) associated with the square node, while the non-diagonal entries depend only on the material property (example: $ k_{m,n}^{m-1,n} $) associated with the edge joining the square node and a circle node. 
Since this pattern repeats over the entire grid, the matrix $\mathbf A $ associated with the two square nodes in Fig. \ref {fig:fdGridNodeMapping}, for example, are identical. Consequently, the solution vectors $ \mathbf {\hat{m}} $  for the two nodes are identical.
With exceptions made at the corner and edge nodes, one can conclude that, for a single material domain, it is sufficient to compute 9 independent columns of $\mathbf M $, independent of the size $N$. A typical set of these 9 independent nodes is illustrated in Fig. \ref{fig:fd9IndependentNodes}. This can be further reduced to 3 independent nodes (one corner, one edge, and one interior node) with appropriate transformations. Exploiting this  \emph {node mapping}, one can compute the SPAI matrix $\mathbf M $  in constant time, independent of  $N$. To the best of our knowledge, this has not been exploited previously to compute SPAI  in classical, or quantum settings. 
If the domain is composed of two materials as illustrated in Fig. \ref{fig:2materials}, then one can show that only 21 independent columns of $\mathbf M $ need to be computed, independent of $N$. This can be further reduced to 10 independent columns through appropriate transformations.
\subsection{Proposed Algorithms}
\label{sec:Algorithms}
To summarize, the proposed algorithm to compute $\mathbf M $  is described in Alg. \ref{alg:SPAI} where
\begin{enumerate}
    \item The sparsity pattern of $\mathbf K$  is first copied over to $\mathbf M$.
    \item Then, for each column of   $\mathbf K$, if that column (node) is mapped to a previously computed column (node), we copy the previously computed solution
    \item Else we call the \emph {sparse box algorithm} (see below), and the computed solution is pushed to $\mathbf M$.
    \item Finally, we force symmetry on $\mathbf M$ to address possible numerical errors.
\end{enumerate}

\begin{algorithm}[H]
    \caption{Sparse Approximate Inverse}
    \label{alg:SPAI}
    \begin{algorithmic}[1]
        \Procedure {QUBOSparseApproximateInverse}{$\mathbf K$}  
        \State  $ \mathbf M \leftarrow \mathbf K  $ \Comment{M has same sparsity as K} \label{alg:center_init}
        \State $N  \leftarrow \text{dim} \mathbf K$ \Comment{number of rows in K}
        \For{$j \leftarrow$ $0$ to $N-1$}  
                \If {nodeMap[j] $<$ j } \Comment{If node is mapped to another node}
                \State $\mathbf M[:,j] \leftarrow  \mathbf M[:,nodeMap[j]] $ \Comment{copy solution}
                \Else 
                \State $\mathbf {\hat m} \leftarrow $ SparseBoxAlg($\mathbf  K$, $j$)   \Comment{solve via QUBO}
                \State $\mathbf s \leftarrow \mathbf K[:,j]$ \Comment{non-zero rows of column j}
                \State $\mathbf M[\mathbf s,j] \leftarrow  \mathbf {\hat m} $ \Comment{copy solution}
                \EndIf
        \EndFor                
        \State $\mathbf M \leftarrow  (\mathbf M + \mathbf M^T)/2 $ \Comment{impose symmetry}
        \EndProcedure    \Comment{Output: solution $ \mathbf M$}
    \end{algorithmic}
\end{algorithm}

The above algorithm uses a sparse box algorithm, which is a generalization of the box algorithm proposed in \cite{srivastava2019box}. The original box algorithm solves small dense linear systems using a QUBO formulation (see \cite{srivastava2019box} for details). It is modified here to solve the sparse problem in Eq. \ref{reducedProblem}. A few observations regarding the proposed sparse box algorithm (see Alg. \ref{alg:Box}) are:
\begin{enumerate}
    \item We have chosen an initial box size of $L = 1$. This is an arbitrary choice; the algorithm is robust for any reasonable value \cite{srivastava2019box}; see numerical experiments. Choosing a large initial value for $L $ will increase the number of contraction steps, while a small initial value will increase the number of translation steps.  
     \item A total of 2s qubits are created, and the initial potential energy is zero.
    \item In the main iteration, the QUBO problem is constructed using the software packages pyQUBO \cite{zaman2021pyqubo}.
    \item The QUBO problem can be solved in exactly or through (3) quantum annealing; see numerical experiments.
    \item If the computed potential energy is less than the current minimum, then we have found a better solution. In this case, we move the box center $\mathbf c$ to the new computed solution (translation) to improve the accuracy. Otherwise, the current solution is optimal, and we shrink the box size $ L $ (contraction) to improve the precision. 
    \item For termination, using a very small value, say, $ \epsilon _{box}  \approx 10^{-14}$, is not desirable since (1) it will increase the computation cost, and (2) as the box size becomes small, the potential energy function $\Pi$ becomes relatively flat, and the computed potential energies from various qubit configurations will be numerically equal; this will result in the algorithm choosing the wrong step (i.e. translate versus contract). Further, $ \epsilon _{box}  \approx 10^{-14}$ is not required since we are only solving for an approximate inverse. A typical choice is $ \epsilon _{box}  \approx 10^{-6}$; see numerical experiments.
\end{enumerate}

\begin{algorithm}[H]
    \caption{Sparse Box Algorithm}
    \label{alg:Box}
    \begin{algorithmic}[1]
        \Procedure {SparseBoxAlg}{$\mathbf K$, $j$}
        \State $\mathbf s \leftarrow \mathbf K[:,j]$ \Comment{non-zero rows of column j}
        \State $s \leftarrow \text{dim} (\mathbf s)$ \Comment{sparsity of column j}
        \State $i \leftarrow  \mathbf (s [i] == j)$ \Comment{find the entry for column j}
        \State  $ \mathbf c \leftarrow \mathbf 0  $ \Comment{center of length $n$} \label{alg:center_init}
        \State  $L \leftarrow 1$ \Comment{initialize box size} \label{alg:length_init}
        \State $\mathbf q_1,\mathbf  q_2   \leftarrow \text{Qubits}(s)   $  \Comment{create qubit arrays of length s}
        \State $ iter = 0$
        \State  $\Pi_{\text{min}} \leftarrow 0  $ \Comment{energy initialization} \label{alg:energy_init}
        \Repeat \Comment{until convergence}
        \State $\mathbf {\hat m} \leftarrow \mathbf c + L(-2\mathbf q_1 +\mathbf q_2)$\Comment{symbolic expression}
        \State $\Pi \leftarrow \frac 1 2 \mathbf {\hat m}^\intercal \mathbf K[\mathbf s,\mathbf s] \mathbf {\hat m}-\mathbf {\hat m}[i] $ \Comment{construct QUBO}
        \State $\Pi^*, \mathbf q_1^*, \mathbf q_2^*  \leftarrow \text{minimize}(H)$ \Comment{solve QUBO}
        \If {$\Pi^*< \Pi_{\text{min}}$}
        \State $ \mathbf c \leftarrow \mathbf c + L(-2\mathbf q_1^* +\mathbf q_2^*) $ \Comment{translation}
        \State $\Pi_{\text{min}} \leftarrow \Pi^*$ \Comment{new minimum}
        \Else 
        \State $ L \leftarrow L/2$ \Comment{contraction} 
        \EndIf
        \State $ iter = iter + 1$
        \Until ({$L< \epsilon _{box} $}) or ({$iter > iter_{max}$}) \Comment{termination}
        \EndProcedure    
        \Comment{Output: solution $\mathbf c$}
    \end{algorithmic}
\end{algorithm}

\subsection{D-WAVE Embedding}
\begin{figure}[H]
     \begin{subfigure}[c]{0.45\textwidth}	
 	\begin{center}
	\includegraphics[width=0.8\textwidth]{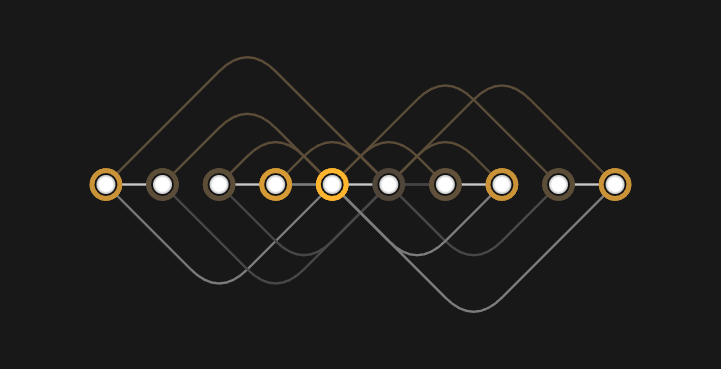}
	\end{center}
        \caption{}
        \label{fig:logicalQubits}
 \end{subfigure}
 \begin{subfigure}[c]{0.36\textwidth}	
        \begin{center}
	\includegraphics[width=1\textwidth]{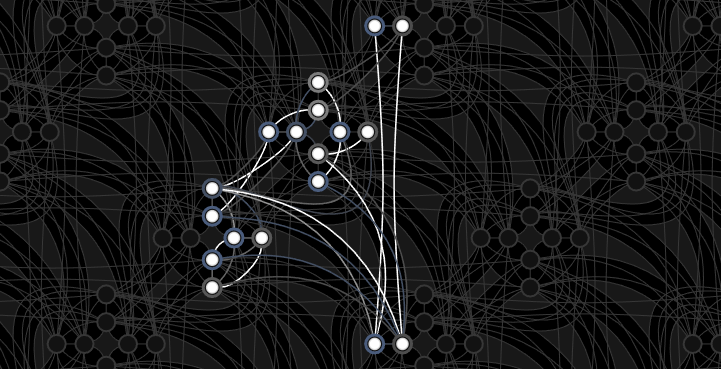}
        \end{center}
        \caption{}
        \label{fig:physicalQubits}
 \end{subfigure}
\caption{(a) 10 logical qubits. (b) 18 physical qubits on Pegasus architecture. }
 \end{figure}

We now consider mapping the QUBO problems generated by  Alg. \ref{alg:Box} to the D-Wave Advantage quantum annealing machine. The D-wave Advantage is equipped with 5000+ qubits, embedded in a Pegasus architecture. Each QUBO problem involves at most 10 logical qubits with the connectivity graph illustrated in Fig. \ref{fig:logicalQubits}. Using the default embedding, these logical qubits are mapped to 18 physical qubits illustrated in  Fig. \ref{fig:physicalQubits}. The default chain strength was found to be sufficient for all numerical experiments.

\section{Numerical Experiments}
\label{sec:Experiments}
In the following experiments, we consider a $(g_x+2) \times (g_y+2)$ rectangular finite-difference grid as shown earlier in Fig. \ref{fig:fdGrid}.  We assume that the field on the boundary is zero, i.e. we only solve for the interior $g_x \times g_y$ grid.  The default values, unless otherwise noted, for all the experiments are:
\begin{itemize}
    \item The grid size is $g_x = 401$ and $g_y = 301$.
    \item The rectangle is composed of a single material with $k = 1$.
    \item The box-tolerance in Algorithm 2 is set to $\epsilon _{box} = 10^{-6}$.
    \item The box length in Algorithm 2 is initialized to $L = 1$.
    \item The maximum box iterations in Algorithm 2 is set to 100.
    \item The conjugate gradient residual tolerance is set to $10^{-10}$
\end{itemize}

In the experiments, Q-PCG refers to the standard preconditioned conjugate gradient method where the proposewd SPAI preconditioner is used. In each of the following experiments, we graph the residual error against the number of iterations for both Q-PCG and CG. The implementation is in Python, and uses pyQUBO \cite{zaman2021pyqubo} to create the QUBO model, and pyAMG \cite {pyamg2023} to construct the $ \mathbf K $ matrix. 

\subsection{Performance Evaluation}

\begin{figure}[H]
     \begin{subfigure}[c]{0.45\textwidth}	
 	\begin{center}
	\includegraphics[width=0.8\textwidth]{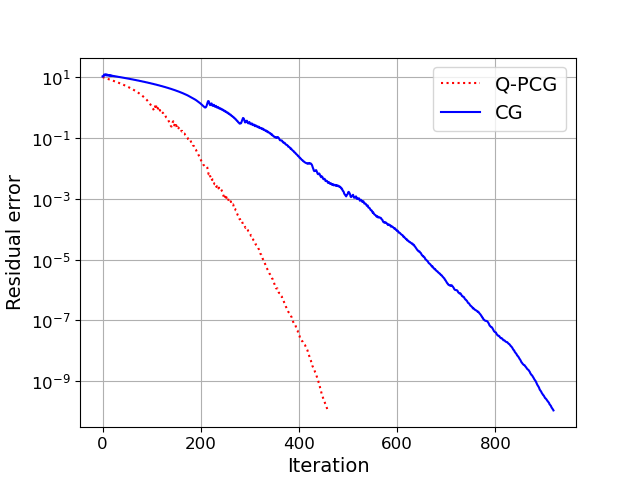}
	\end{center}
        \caption{}
        \label{fig:expt1a}
 \end{subfigure}
 \begin{subfigure}[c]{0.45\textwidth}	
        \begin{center}
	\includegraphics[width=1\textwidth]{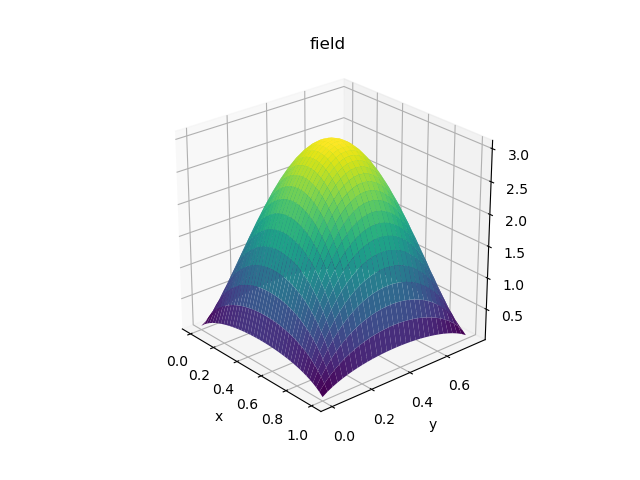}
        \end{center}
        \caption{}
        \label{fig:expt1b}
 \end{subfigure}
\caption{(a) Convergence plot for $g_x = 401$, $g_y = 301$. (b) Computed Poisson field. }
 \end{figure}

In the first experiment, we compare the convergence of Q-PCG and CG using the default values listed above; this results in $N = 120,701$ (size of the $ \mathbf K$  matrix).  The convergence plots using  D-WAVE's \emph {dimod} exact QUBO solver are illustrated in Fig. \ref{fig:expt1a}. \emph{Exactly the same Q-PCG convergence plot was obtained when using the D-WAVE quantum annealing solver.} The SPAI does not significantly improve the convergence of CG here since $\mathbf K$ is relatively well-conditioned. The computed Poisson field is illustrated in Fig. \ref{fig:expt1b}.

The two solvers are compared in Table \ref{tab:tableTiming}.  

\begin{table}[h!]
    \captionsetup{width=1\textwidth}
    \caption{Comparing the two solvers.} 
    \label{tab:tableTiming}
    \begin{tabular}{l |c|r} 
      \textbf{   } & \textbf{Exact} & \textbf{Quantum}\\
      \hline
      Total QUBO solves &  294 & 294  \\
      Avg. time per QUBO solve &  0.9 milliseconds & 26 milliseconds \\
      Time to compute $ \mathbf M$ & 0.5 seconds &  209 seconds \\
    \end{tabular}
\end{table}
Additional details on the QPU timing for a single QUBO solve (using the default 100 samples) are provided in Table \ref {tab:tableQPUTiming}. As one can observe there is significant overhead in D-WAVE QPU allocation, programming, access and post-process. 
\begin{table}[h!]
    \captionsetup{width=1\textwidth}
    \caption{QPU timing per QUBO solve.} 
    \label{tab:tableQPUTiming}
    \begin{tabular}{l |r} 
      \textbf{  \textbf{Task} } &  \textbf{Milliseconds}\\
      \hline
      Access &  26 \\
      Programming &  15 \\
      Sampling  &  11  \\
      Readout &  7  \\
      Post-process & 2 \\
      Anneal  &  2 \\
      Delay  &  2 \\
      Access overhead & 1.6 \\
    \end{tabular}
\end{table}

\subsection{Multiple Materials}
\begin{figure}[H]
   \begin{subfigure}[c]{0.45\textwidth}	
 	\begin{center}
	\includegraphics[width=0.8\textwidth]{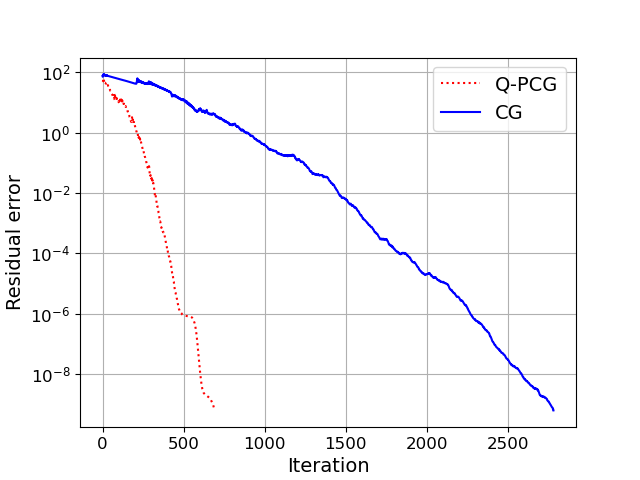}
	\end{center}
        \caption{}
        \label{fig:expt3a}
 \end{subfigure}
 \begin{subfigure}[c]{0.45\textwidth}	
        \begin{center}
	\includegraphics[width=0.8\textwidth]{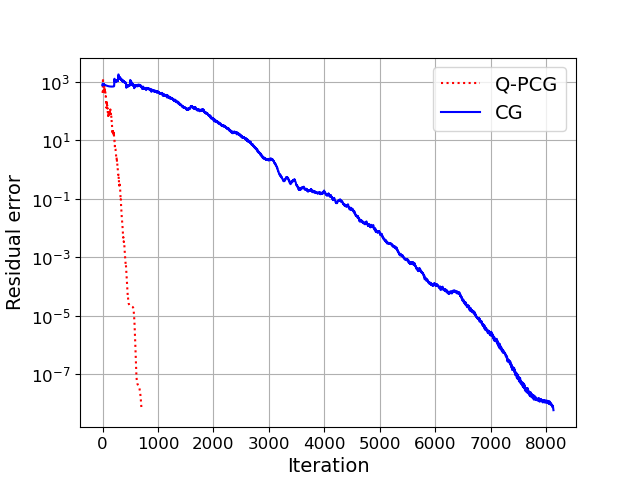}
        \end{center}
        \caption{}
        \label{fig:expt3b}
 \end{subfigure}
  \caption{Convergence plots for two-material domain: (a)  $ k_1 = 1, k_2 = 10$. (b) $k_1 = 1, k_2 = 100$. }
 \end{figure}

The real advantage of the SPAI becomes evident when we consider two materials with the $k_1$ on the left and $k_2$ on the right (see Fig. \ref{fig:2materials}). For the default grid size of $g_x = 401$ and $g_y = 301$, we observe in Fig. \ref{fig:expt3a} that  Q-PCG provides a speed-up of 3.3 when $k_1 = 1$ and $k_2 = 10$, and a speed-up of 11 when $k_1 = 1$ and $k_2 = 100$. Such multi-material problems are fairly common in engineering, especially during topology optimization \cite{zuo2017multi, suresh2013efficient}.

\subsection{Impact of box tolerance}
We now study the impact of the box tolerance on the convergence of Q-PCG. For a single material problem, with default values, Table \ref{tab:table1} summarizes the Q-PCG iterations and total box iterations using the exact and quantum annealing solvers (recall from Fig. \ref{fig:expt1a} that regular CG converges in about 900 iterations). As the box tolerance is varied between $10^{-8}$ and $10^{-2}$, one can observe in Table \ref {tab:table1} that the number of Q-PCG iterations remains around 459, while the number of box iterations decreases, as expected. However, Q-PCG did not converge to the desired tolerance, i.e., the SPAI matrix is not effective when the box tolerance is too coarse ($10^{-1}$).

\begin{table}[h!]
    \captionsetup{width=1\textwidth}
    \caption{Impact of box tolerance on Q-PCG iterations and (total box iterations).} 
    \label{tab:table1}
    \begin{tabular}{l |c|r} 
      \textbf{$\epsilon _{box} $  } & \textbf{Exact} & \textbf{Quantum}\\
      \hline
      $10^{-8}$ & 458 (401) & 459 (401)  \\
      $10^{-6}$ & 459 (294) & 459 (294) \\
      $10^{-4}$ & 459 (202) & 459 (202) \\
      $10^{-2}$ & 464 (96) & 473 (96) \\
      $10^{-1}$ & - (50) & - (50)  \\
    \end{tabular}
\end{table}

\subsection{Impact of box length}
Next we study the impact of the box length on the convergence of Q-PCG. For a single material problem, with default values, Table \ref{tab:table2} summarizes the Q-PCG iterations and total box iterations using exact and quantum annealing solvers. As the box length is varied between $10^{-2}$ and $10^{4}$, one can observe in Table \ref {tab:table2} that the number of Q-PCG iterations remains at 459, while the number of box iterations is minium when $L = 1$. A star (*) indicates that the desired box tolerance was not achieved, and the box algorithm exited when the maximum iteration was reached. However, even in this case, Q-PCG converged.
\begin{table}[h!]
    \captionsetup{width=1\textwidth}
    \caption{Impact of box length on Q-PCG iterations and (total box iterations). }
    \label{tab:table2}
    \begin{tabular}{l |c|r} 
      \textbf{ $L $} & \textbf{Exact} & \textbf{Quantum}\\
      \hline
      $10^{4}$ & 458 (437) & 459 (437)  \\
      $10^{2}$ & 458 (390) & 459 (365) \\
      $10^{1}$ & 458 (315) & 458 (315) \\
      $10^{0}$ & 458 (294) & 458 (294) \\
      $10^{-1}$ & 458 (333) & 458 (333) \\
      $10^{-2}$ & 691 (431*) & 691 (431*)  \\
    \end{tabular}
\end{table}

\section{Conclusions}
A hybrid classical-quantum strategy for solving large linear systems of equations was proposed. The strategy relied on computing a sparse approximate preconditioner (SPAI) on a quantum annealing machine and using this preconditioner, along with an iterative solver, on a classical machine. Its effectiveness was demonstrated on large ($ N > 100,000$) ill-conditioned matrices arising from a finite-difference formulation of the Poisson problem. 

There are many directions for continued research. (1) While we demonstrated the method using quantum annealing machines, it can be extended to quantum gate computers since QUBO problems can be solved (approximately) on such machines via quantum approximate optimization algorithms \cite{farhi2014quantum, clader2013preconditioned}. (2) The finite difference formulation of the Poisson problem resulted in matrices with small sparsity ($s = 5$). Other formulations and other field problems would lead to a larger $s$. For example, 3D structured-grid finite element analysis of elasticity problems would result in $s = 27 \times 3 = 81$. The proposed strategy and algorithm, in theory, extends to such scenarios. However, the performance of SPAI and efficient embedding of the QUBO problems need to be investigated. (3) We selected a simple \emph {a priori} sparsity pattern for the preconditioner $M$; adaptive patterns and their impact on Q-PCG need to be explored. (4) One can potentially exploit parallel annealing \cite{pelofske2022parallel} for computing $M$. (5) We limited the box representation to 2 qubits per real variable (see Eq. \ref{boxRepresentation}); extending to multiple qubits will accelerate convergence but will require a larger number of qubits. 

\section*{Compliance with ethical standards}
\label{sec:ethics}
The authors declare that they have no conflict of interest.

\section*{Replication of Results}
\label{sec:replic}
The Python code pertinent to this paper is available at \href{https://github.com/UW-ERSL/SPAI/}{https://github.com/UW-ERSL/SPAI}.

\section*{Acknowledgments}
We would like to thank the Graduate School of UW-Madison for the the Vilas Associate grant.

\bibliography{QCLinearSolver,ERSLReferences}

\end{document}